\documentclass[12pt,a4paper,reqno,oneside]{amsart}
\usepackage{config}

\usepackage[style=alphabetic,backend=biber,maxbibnames=6]{biblatex}
\NewBibliographyString{toappear}
\DefineBibliographyStrings{english}{%
  toappear = {to appear},
}

\renewbibmacro*{in:}{%
  \iffieldundef{pubstate}
    {}
    {\printfield{pubstate}%
     \setunit{\addspace}%
     \clearfield{pubstate}}%
  \bibstring{in}%
  \printunit{\intitlepunct}}

\graphicspath{{./pics/}}
\begin{document}
\title{Frames and representing systems in Hardy spaces}
\author{Timur Batenev}
\address{Saint Petersburg University, 7/9 Universitetskaya nab., St. Petersburg, 199034, \mbox{Russia}}
\email{tbatenev00@gmail.com}


\begin{abstract}
 We study representing systems of Cauchy kernels or dilations of a given function
 in Hardy spaces. We give a complete geometric description
 of countable sets $\Lambda \subset \D$ such that the system $\{\hat{k}_{\lambda}\}_{\lambda \in \Lambda}$
 of normalized Cauchy kernels forms a generalized frame
 in the Hardy space $H^p(\mathbb{D})$ with $p \in (1, \infty)$ with respect to a particular coefficient space.
 For a wide class of functions we construct representing systems of their dilations in $H^2$.
 We also obtain a linear method for expanding functions in series of Cauchy kernels in $H^2$.     
\end{abstract}

\maketitle

\section{Introduction}

The notion of a basis in a finite-dimensional vector space $X$ is well studied and fully understood. A basis is a 
linearly independent system $X' \subset X$ that spans $X$, i.e., $\text{span } X'=X$. On the other hand, in 
the infinite-dimensional separable case, there are several non-equivalent notions
 of a basis even for a Hilbert space 
(orthonormal basis, Riesz basis, Schauder basis...). In the infinite-dimensional case, completeness
 of a system $X'$ (i.e., $\text{Clos span } X'=X$) does not imply that every element of $X$ can be represented by a convergent
series of scalar multiples of elements of $X'$. 
If such an expansion does exist for every element of $X$, then the system $X'$
 is called a representing system.

\begin{definition}
    A system of elements $\{f_n\}$ of an infinite-dimensional separable complex topological vector space $F$ is called
     a \textit{representing
    system} for $F$ if for any $f \in F$ there exists a sequence $\{c_n\} \subset \C$ such that 
    $$f=\sum\limits_{n}c_n f_n$$
    and the series converges in the topology of $F$.
    If, in addition, the space $F$ is a normed space and the series converges absolutely,
     then the system $\{f_n\}$ is called an \textit{absolutely representing system}.   
\end{definition}

This definition is very similar to that of a Schauder basis but does not require the uniqueness of the coefficients. However, unlike Schauder bases, representing systems (even absolutely representing systems)
 exist in every separable Banach space.
The first systematic theory of representing systems was developed by Yu.\,F. Korobeinik and his school (for Fr\'echet spaces), see \cite{Kor81, Kor87}. It grew out of earlier studies of concrete systems; in particular,
 A.\,F. Leont'ev studied the representation of analytic functions by series of exponentials
(in spaces of analytic functions on domains in $\C$); see \cite{Leo69, Leo72}.

It is known that there are representing systems in every separable Fr\'echet space \cite{KK92, BFGR17}. In Section 2 we give 
an explicit method for constructing representing systems in separable quasi-normed spaces. 

Let $H$ be a reproducing kernel Hilbert space (RKHS) on a set $M$ with reproducing kernel $k_\mu$, and let
$\Lambda \subset M$.
The relationship between the geometric properties of the system $\{k_{\lambda}\}_{\lambda \in \Lambda}$ and the set $\Lambda$
is of great interest and importance. There is a vast literature on this topic for 
RKHS of analytic functions.  

It is well known that in many classical spaces of analytic functions
 (the Paley--Wiener space, the Fock space, the Bergman space)
 sequences of normalized reproducing kernels can form frames in the sense
of Duffin--Schaeffer or even a Riesz basis. For instance, a complete description of such frames
 in the Bergman space was given 
by K. Seip \cite{Sei93}; for further results, see \cite{BDK07, Sei13}. However, there is no complete and
simultaneously minimal system of reproducing kernels (Cauchy kernels)
in the Hardy space $H^2=H^2(\D)$. Hence, there are no Riesz bases of reproducing kernels.
 Moreover, there is no 
Duffin--Schaeffer frame of reproducing kernels.
So it is natural to ask whether they can form a representing system.
In 2018, it was proved by E. Fricain, L.\,H. Khoi and P. Lef\`evre \cite{FKL18}
 that there is no absolutely representing system of reproducing kernels
for $H^2$ and many other classical spaces. The first examples of representing systems of Cauchy
kernels for $H^2$ were constructed by K. S. Speransky and P. A. Terekhin in 2018 \cite{ST18}. 
Their proofs are based on abstract methods of functional analysis. Later A. Baranov
 and the author presented an elementary construction of representing systems of Cauchy kernels for $H^p$, $1\le p<\infty$, and for the disc-algebra, see \cite{BB23}. 

The notion of a basis in a Hilbert space is closely related to the notion of a frame.
 There are also several non-equivalent
notions of a frame in a Banach space. We use the following:

\begin{definition}
Let $F$ be a Banach space, let $\mathcal{F}=\{f_n\}_{n}\subset F$, and let $X$
 be a Banach space of sequences such that
 the unit vectors $e_k(n)=\delta_{kn}$
form a Schauder basis for $X$ (hence $X^*$ is isomorphic to a sequence space). 
We say that a triple $(F, \mathcal{F}, X)$ 
is a \textit{frame} if there exist constants $A, B>0$ such that $\{g(f_n)\}_n \in X^*$
for every $g \in F^*$ and the following frame inequality holds:

$$ A\|g\|_{F^*} \le \|\{g(f_n)\}\|_{X^*}\le B \|g\|_{F^*}.$$
\end{definition}

This notion of a generalized frame was introduced by P.\,A. Terekhin in \cite{Ter09}, 
where it was also shown that every
frame is a representing system. Moreover, in this case, the expansion coefficients may be chosen
 from the space $X$. The systems constructed in \cite{ST18, ST19}
were frames of normalized Cauchy kernels
associated with the space $X=\Big(\bigoplus\limits_{k=1}^{\infty}\ell^2_{n_k}\Big)_{\ell^1}$. Namely, let $n_k \in \N$ and $r_k\in (0,1)$ with $r_k \to 1-$ as $k \to \infty$, and define the sequence 
$\Lambda$ by
\begin{equation}
    \Lambda = \{\lambda_{k,j}=r_ke^{2\pi i j/n_k}: k\ge  1, j=1,...,n_k\}.
    \label{eq:lambda-grid}
\end{equation}
If there exist positive constants $A$ and $B$ such that $A \le n_k(1 - r_k) \le B$
for all $k$, then the system $\{{k}_{\lambda_{k,j}}/\|{k}_{\lambda_{k,j}}\|_2\}_{k=1, j=1}^{\infty, n_k}$ is a frame for $H^2$ with respect to $X$.
The sequence $\Lambda$ is assumed to have the standard
lexicographical order: $\lambda_{1,1}, \,\lambda_{2,1},\, \lambda_{2,2}, \, \lambda_{3,1} \ldots$ . In what follows we always assume that
sequences with double (or triple) indices are ordered in this way.

Let $\Lambda=\{\lambda_n\}_{n=1}^{\infty}\subset \D$, $k_{\lambda}(z)=(1-\overline{\lambda} z)^{-1}$
be the Cauchy kernel at the point $\lambda$ and $\hat{k}_{\lambda}(z) =
 k_{\lambda}(z)/\|k_{\lambda}\|_{H^p}$ be the normalized kernel. Set 
$\mathcal{K}_p(\Lambda)=\{\hat{k}_{\lambda}\}_{\lambda \in \Lambda}$. The first question studied in the paper is the following: 

\begin{Question}
     How can one geometrically describe the sets $\Lambda$ for which the corresponding system of normalized Cauchy kernels $\mathcal{K}_p(\Lambda)$ is a frame for $H^p$ associated with
     $\Big(\bigoplus\limits_{k=1}^{\infty}\ell^p_{n_k}\Big)_{\ell^1}$ with respect to some partition $\Lambda= \bigcup\limits_{k} \Lambda_k$?
\end{Question}
 We give an answer to this question in Section 3; namely, we prove the following

 \begin{Theorem} 
  Let $p\in (1,\infty)$, $\Lambda=\{\lambda_n\}_n \subset \D$. Then there exists
   a partition $\Lambda=\bigcup\limits_{k}\Lambda_{k}$ such that $n_k=\text{card }\Lambda_k < \infty$ and the corresponding system
  $\mathcal{K}_p(\Lambda)$ is a frame for $H^p$ with respect to the space
  $\Big(\bigoplus\limits_{k}\ell^p_{n_k} \Big)_{\ell^1}$ if and only if there exist constants $C,c>0$
   such that for every $r \in (0,1)$ and every nonempty angle
    $$\Gamma=\{z\ \in \C: \text{arg }z \in (\alpha, \beta)\}$$
    there exists $\Lambda'\subset \Lambda \cap \Gamma \cap \{r<|z|<1\}$ with the following properties: 

    (1) the Carleson constant of $\Lambda'$ is bounded by $C$;
     
    (2) $\sum\limits_{\Lambda'}(1-|\lambda|)\ge c|\Gamma|$, where $|\Gamma|$ is the length of the arc $\Gamma \cap \T$.  
 \end{Theorem}
 
 Representing systems of Cauchy kernels are of the form $\{f(\overline{\lambda}_n z)\}_{n}$; that is,
 they are dilations of the function $f(z)=(1-z)^{-1}$.  
 In Section 4, using the discretization technique introduced by A. Baranov and the author in \cite{BB23}, we prove the following
\begin{Theorem}
     Suppose that a function $g(z)=\sum\limits_{n \ge 0} \beta_n z^n$ is analytic in the unit disk, $\beta_n \neq 0$ and $\lim\limits_{n \to \infty}|\beta_n|^{1/n}=1.$ Then there exists a set $\{\lambda_n\} \subset \mathbb{D}$ such that $\{g(\overline {\lambda}_n z)\}$ is a representing system for $H^2$.
\end{Theorem}
For instance, the condition holds for the following functions:
 
 1) $g_1(z)=(1-z)^{-\alpha}$, $\alpha \in \R \setminus\{0, -1, -2,...\}$; then $g_1(\bar\lambda z)$ is the reproducing kernel for the weighted Bergman space $A^2_{\alpha-2}$ ($\alpha>1$);

 2) $g_2(z)=-\log(1-z)/z$; then $g_2(\bar\lambda z)$ is the reproducing kernel for the Dirichlet space; 

 3) $\displaystyle g_3(z)=\exp\bigg(a\frac{1+z}{1-z}\bigg)$, $a>0$,
 which is the reciprocal of a singular inner function.

 Indeed, for $g_3(z)=\sum\limits_{n\ge 0}\beta_n z^n$ one has $g'_3(z)/g_3(z)=2a/(1-z)^2$;
 hence $(1-2z+z^2)g'_3(z)=2ag_3(z)$. Therefore, 
 $$\beta_{n+1}=2\frac{n+a}{n+1}\beta_n-\frac{n-1}{n+1}\beta_{n-1}, \quad n\ge 1,$$
 and $\beta_0=e^a,\, \beta_1=2ae^a$, which yields $\beta_{n+1} \ge \beta_n \ge \beta_1 \ge 2ae^a$ for $n\ge 1$,
 and $\liminf\limits_{n} \beta_n^{1/n}\ge 1$.
 On the other hand, $\limsup\limits_{n} \beta_n^{1/n}\le 1$ as $g_3$ is analytic in the unit disk. 
  
 Suppose that $\{f_n\}_n$ is a representing system for a Banach space $F$. Another interesting question
 is whether it is possible to find linear functionals $\{c_n(\cdot)\}_{n} \subset F^*$ 
 such that $f=\sum\limits_{n}c_n(f)f_n$. In this case we say that the representing system \textit{allows a linear method for
the choice of coefficients}. 

In Section 5 we prove that the frames constructed by K. S. Speransky and P. A. Terekhin do not allow a
linear method (in a somewhat stricter sense) and prove the following 
 \begin{Theorem}
      Define $\Lambda$ as in \eqref{eq:lambda-grid} and let $n_k(1-r_k) \to \infty$.
    Then the system $\mathcal{K}(\Lambda)$ is a representing system for
    $H^2$ and it allows a linear method for the choice of coefficients.
 \end{Theorem}

\section{A general construction of representing systems in quasi-normed spaces}
\begin{Proposition}
     Let $X$ be a separable quasi-normed space with a constant $C\ge 1$;
     that is $$\|x+y\| \le C(\|x\|+\|y\|) \text{ for all } x, \,y \in X.$$
     Let $\{x_n\}_{n=1}^{\infty}$ be a complete system in $X$, and let a sequence
$\{y_j\}_{j=1}^{\infty}$ satisfy the following condition: $y_j\neq y_k$ for any $j\neq k$, and $\text{for every } n \text{ there exists a subsequence of }
\{y_j\}$ that converges to $x_n$. Then $\{y_j\}_j$ is a representing system for $X$. 
\end{Proposition}

\begin{proof}
Fix an enumeration $\{y_j\}$. We write $x<y$ if there exist $j,k: x=y_j, y=y_k$ and $j<k$.
    
Let $z=z_0 \in X$, $M_0=0$. Choose a linear combination $\sum\limits_{n=1}^{N_0}c_{0,n}x_n$ such that 
$\Big\|\sum\limits_{n=1}^{N_0}c_{0,n}x_n-z_0 \Big\| < 1/(2C).$ 
For some $\varepsilon_1$ (to be specified later) 
we can find $\{y_{0,n}^0\}_{n=1}^{N_0}$ such that $y_{0,1}^0<y_{0,2}^0<...<y_{0,N_0}^0$
and for every $n=1,...,N_0$ one has $\|x_n-y_{0,n}^0\|\leq \varepsilon_1$. 
Then we get $$\bigg \|\sum\limits_{n=1}^{N_0}c_{0,n}y_{0,n}^0-z_0 \bigg\| \leq 
C \bigg \|\sum\limits_{n=1}^{N_0}c_{0,n}(x_n-y_{0,n}^0)\bigg \|+C \bigg \|\sum\limits_{n=1}^{N_0}c_{0,n}x_n-z_0 \bigg\|$$ 
$$<C\varepsilon_1 \bigg (\sum\limits_{n=1}^{N_0-1}C^n|c_{0,n}|+C^{N_0-1}|c_{0,N_0}|\bigg)
+C \bigg \|\sum\limits_{n=1}^{N_0}c_{0,n}x_n-z_0 \bigg\| <1/2$$ 
by the choice of $\varepsilon_1$.

Now define $z_1=z_0-\sum\limits_{n=1}^{N_0}c_{0,n}y_{0,n}^0$.
 Proceeding in this way we obtain $z_2,...,z_k$. Choose the linear combination $\sum\limits_{n=1}^{N_{k}}c_{k,n}x_n$ such that 
$\|\sum\limits_{n=1}^{N_k}c_{k,n}x_n-z_k\| < 1/(2^{k+1}C).$ Put 
 $$M_k=2^{k+1}([C(\sum\limits_{n=1}^{N_k-1}C^n \|x_n\| |c_{k,n}|+
C^{N_k-1}\|x_{N_k}\||c_{k,N_k}|)]+1),$$
where $[\cdot]$  denotes the floor function. 
For some $\varepsilon_k$ (to be specified later) we can find 
$Y^{m}_k=\{y_{k,n}^m\}_{n=1}^{N_{k}}, m=1,...,M_k$ such that for every $n=1, ..., N_{k}$ and  
$m=1,..., M_k$ the inequality  $\|x_n-y_{k,n}^m\|\leq \varepsilon_k$ holds and
the elements are ordered in the following way: within each group $Y_k^{m}$ the elements are ordered
 naturally, while $y_{k,1}^1>y_{k-1,N_{k-1}}^{M_{k-1}}$ and
 $y^{m_1}_{k,n_1}<y^{m_2}_{k,n_2}$ if  $m_1<m_2$.

Then we have 
\begin{align*}
    \bigg\|\frac{1}{M_k}\sum\limits_{m=1}^{M_k}\sum\limits_{n=1}^{N_{k}}c_{k,n}y_{k,n}^m-z_k \bigg\|
    \leq C \bigg \|\frac{1}{M_k}\sum\limits_{m=1}^{M_k}\sum\limits_{n=1}^{N_{k}}c_{k,n}y_{k,n}^m
    -\sum\limits_{n=1}^{N_k}c_{k,n}x_n \bigg\| \\
+C\Big\|z_k-\sum\limits_{n=1}^{N_k}c_{k,n}x_n \Big\| 
\leq \frac{C}{M_k}\bigg (\sum\limits_{m=1}^{M_k-1}C^m 
\bigg\|\sum\limits_{n=1}^{N_k} c_{k,n}(y_{k,n}^m-x_n)\bigg\| \\
+C^{M_k-1}\Big\|\sum\limits_{n=1}^{N_k} c_{k,n}(y_{k,n}^{M_k}-x_n)\Big\|\bigg)
+C\Big\|z_k-\sum\limits_{n=1}^{N_k}c_{k,n}x_n \Big\|=A.
\end{align*}

For every $m=1,..., M_k$ we have $$ \bigg\|\sum\limits_{n=1}^{N_k} c_{k,n}(y_{k,n}^{m}-x_n)\bigg\|
\leq \varepsilon_k \bigg(\sum\limits_{n=1}^{N_k-1}C^n |c_{k,n}|+C^{N_k-1}|c_{k,N_k}|\bigg).$$ Therefore, 
\begin{align*}
A \leq \frac{\varepsilon_k C}{M_k}
\bigg(\sum\limits_{n=1}^{N_k-1}C^n |c_{k,n}|+&C^{N_k-1}|c_{k,N_k}|\bigg)
\bigg(\sum\limits_{m=1}^{M_k-1}C^m+C^{M_k-1}\bigg) \\ 
+&C \bigg\|z_k-\sum\limits_{n=1}^{N_k}c_{k,n}x_n \bigg\| <1/2^{k+1}
\end{align*}

by the choice of $\varepsilon_k$. Moreover, for $1 \le M \le M_k$ and $1 \le N \le N_k$, we obtain
$$\bigg \|\frac{1}{M_k}\sum\limits_{m=1}^{M-1}\sum\limits_{n=1}^{N_k}c_{k,n}y_{k,n}^m+
\frac{1}{M_k}\sum\limits_{n=1}^{N}c_{k,n}y_{k,n}^M-z_k \bigg \|\leq \frac{C}{M_k}
\bigg \|\sum\limits_{m=1}^{M-1}\sum\limits_{n=1}^{N_k}c_{k,n}y_{k,n}^m-(M-1)z_k \bigg\|$$
$$+C^2\bigg \|\frac{M_k-(M-1)}{M_k}z_k\bigg \|
+\frac{C^2}{M_k}\bigg \|\sum\limits_{n=1}^{N}c_{k,n}y_{k,n}^M\bigg\|=B.$$
 One can estimate the first term in $B$ as above:
\begin{align*}
\frac{C}{M_k} \bigg\|\sum\limits_{m=1}^{M-1}\sum\limits_{n=1}^{N_k}&c_{k,n}y_{k,n}^m-(M-1)z_k \bigg \| \\ 
\leq \frac{\varepsilon_k C}{M_k}
\bigg(\sum\limits_{n=1}^{N_k-1}C^n |c_{k,n}|+&C^{N_k-1}|c_{k,N_k}|\bigg)
\bigg(\sum\limits_{m=1}^{M_k-1}C^m+C^{M_k-1}\bigg).
\end{align*}
For the second term we get
$C^2\Big\|\frac{M_k-(M-1)}{M_k}z_k \Big\|\leq C^2/2^k$, and for the last term in $B$ we have:
$$\bigg\|\sum\limits_{n=1}^{N}c_{k,n}y_{k,n}^M \bigg\| 
\leq C \bigg\|\sum\limits_{n=1}^{N}c_{k,n}(y_{k,n}^M-x_n)\bigg \|
+C\bigg\|\sum\limits_{n=1}^{N}c_{k,n}x_{n}\bigg\|$$
$$\le \varepsilon_k C \bigg(\sum\limits_{n=1}^{N-1}C^n |c_{k,n}|+C^{N-1}|c_{k,N}|\bigg)+
C\bigg(\sum\limits_{n=1}^{N-1}C^n \|x_n\| |c_{k,n}|+C^{N-1}\|x_N\||c_{k,N}|\bigg).$$
Finally, $B \leq 4C^2/2^{k+1}$ 
by the choice of $\varepsilon_k$. Now we define $z_{k+1}=z_k-\frac{1}{M_k}\sum\limits_{m=1}^{M_k}\sum\limits_{n=1}^{N_k}c_{k,n}y_{k,n}^m$.

We have shown that the series converges to $z$.
\end{proof}

\begin{Corollary}
Let $\Omega$  be a domain in $\R^d$, and let $H$ be an RKHS consisting of continuous functions on $\Omega$ with the kernel $K(x,y)=k_y(x) \in C(\Omega \times \Omega)$. Suppose that $k_{y_1}\neq k_{y_2}$ whenever $y_1, y_2 \in \Omega$ and $y_1 \neq y_2$.
Then there exists a representing system of reproducing kernels for $H$.       
\end{Corollary}  

\begin{Corollary}
    Let $p \in (0, \infty)$, and let $\{\mu_n\}_n$ be a non-Blaschke sequence in $\D$,
    $\lambda_{n,j} \in \D, \lambda_{n_1,j_1} \neq \lambda_{n_2,j_2}$ if $(n_1, j_1)\neq (n_2,j_2)$, and $\lambda_{n,j} \to \mu_n$ as $j \to \infty$.
Then $\{k_{\lambda_{n,j}}(z)\}$ \textup{(}with any ordering\textup{)} is a representing system for $H^p(\D)$.
\end{Corollary}  

\begin{Corollary}
    If $p\in (0,\infty)$ and $\Lambda=(\mathbb{Q} + i \mathbb{Q})\cap \{|z|<r\}$, then $\{k_{\lambda}(z)\}_{\lambda \in\Lambda}$ is a representing system for $H^p(\D)$ for any $r\in(0,1)$.
\end{Corollary} 
\section{\texorpdfstring
  {Frames of Cauchy kernels in $H^p$ associated with $\Big(\bigoplus\limits_{k=1}^{\infty}\ell^p_{n_k}\Big)_{\ell^1}$-coefficients}
  {Frames of Cauchy kernels in Hardy spaces}}
Recall that $$X_p=\Big(\bigoplus\limits_{k=1}^{\infty}\ell^p_{n_k}\Big)_{\ell^1}= \Bigg\{c=\{c_{k,j}\}^{\infty, n_k}_{k,j=1}:\|c\|= \sum\limits_{k=1}^{\infty}\Bigg(\sum\limits_{j=1}^{n_k}|c_{k,j}|^p\Bigg)^{1/p}< \infty\Bigg\}$$
is a Banach space and the standard unit vectors form a Schauder basis
 for $X_p$, $p \in (1, \infty)$. Therefore, for $1/p+1/q=1$, 
$$X^*_p=\Big(\bigoplus\limits_{k=1}^{\infty}\ell^q_{n_k}\Big)_{\ell^\infty}= \Bigg\{c=\{c_{k,j}\}^{\infty, n_k}_{k,j=1}:\|c\|= \sup\limits_{k\ge 1}\Bigg(\sum\limits_{j=1}^{n_k}|c_{k,j}|^q\Bigg)^{1/q}< \infty\Bigg\}.$$

We use the following notation:

 1) If a system of normalized Cauchy kernels $\mathcal{K}_p(\Lambda)$ is a frame associated with
$\Big(\bigoplus\limits_{k=1}^{\infty}\ell^p_{n_k}\Big)_{\ell^1}$, then $\Lambda_k$ is the subset of $\Lambda$
corresponding to the $k$-th direct summand of $\Big(\bigoplus\limits_{k=1}^{\infty}\ell^p_{n_k}\Big)_{\ell^1}$. Notice that the partition need not be unique;

2) $A_{r,r'}=\{z\in \C:r\le|z|<r'\}$;

3) The Carleson square is
\[
 S(\zeta,h)=\{z\in\D:1-h<|z|<1,\ |\text{arg}(z\bar\zeta)|<h/2\},
 \quad \zeta\in\T,\quad 0<h<1.
\]
For a discrete set $\Lambda \subset\D$, let
$\mu_{\Lambda}=\sum\limits_{\lambda\in\Lambda}(1-|\lambda|)\delta_\lambda$ and
\[
 \text{Carl}(\Lambda) =\sup_{\zeta,h}\frac{\mu_{\Lambda}(S(\zeta,h))}{h}.
\]

Let $\Lambda=\bigcup\limits_{k}\Lambda_k$, $\Lambda_k=\{\lambda_{k,j}\}_{j=1}^{n_k}$. According to the definition, $(H^p, \mathcal{K}_p (\Lambda), X_p)$ is a frame if there exist $A, B > 0$ such that 
$$ A\|\Phi\|_{(H^p)^*}\le \sup\limits_{k\ge 1} \bigg(\sum\limits_{j=1}^{n_k}|\Phi(\hat{k}_{\lambda_{k,j}})|^q\bigg)^{1/q} \le B\|\Phi\|_{(H^p)^*}$$ 
for any $\Phi \in (H^p)^*$.

Recall that for $\Phi \in (H^p)^*$ there exists a unique $g \in H^q$ such that $\Phi(f)=\int\limits_{\T}f \bar{g}dm$ and $\|\Phi\|_{(H^p)^*}\asymp \|g\|_q$. Note that
 $\|k_\lambda\|_p \asymp  (1-|\lambda|)^{-1/q}$ and for any $g \in H^q$ one has $\int\limits_{\T}g \bar{k}_{\lambda}dm=g(\lambda)$. Therefore, to verify that 
 $(H^p, \mathcal{K}_p (\Lambda), X_p)$ is a frame we need to show that 
 $$A\|g\|_{q}\le \sup\limits_{k\ge 1} \bigg(\sum\limits_{j=1}^{n_k}(1-|\lambda_{k,j}|)|g(\lambda_{k,j})|^q\bigg)^{1/q} \le B\|g\|_{q}$$
 for some $A,B>0$ and any $g \in H^q$.

 \begin{Remark}
    If the system $\mathcal{K}_p(\Lambda)$ is a 
    $\Big(\bigoplus\limits_{k=1}^{\infty}\ell^p_{n_k}\Big)_{\ell^1}$-frame
    with respect to a partition $\{\Lambda_k\}_{k}$, then the Carleson constants $\textup{\text{Carl}}(\Lambda_k)$ are uniformly bounded.
    Obviously, this is a necessary and sufficient condition for the upper frame inequality.
\end{Remark}

\begin{Remark}
    If for some partition $\Lambda=\bigcup\limits_{k} \Lambda_k$ the system $\mathcal{K}_p(\Lambda)$ is a frame associated with
    $\Big(\bigoplus\limits_{k=1}^{\infty}\ell^p_{n_k}\Big)_{\ell^1}$ and $ \tilde{\Lambda}=\{\tilde\lambda_n\}_n$
    is an arbitrary countable subset of
    $\D$, then there is a partition of $\Lambda \cup \tilde\Lambda$ such that $\mathcal{K}_p(\Lambda \cup \tilde\Lambda)$
    is a frame associated with $\Big(\bigoplus\limits_{k=1}^{\infty}\ell^p_{m_k}\Big)_{\ell^1}$ for some $m_k$. 
    Indeed, if the frame inequality holds for $\mathcal{K}_p(\Lambda)$ and $\Big(\bigoplus\limits_{k=1}^{\infty}\ell^p_{n_k}\Big)_{\ell^1}$, then it is obviously fulfilled for
    $\mathcal{K}_p(\Lambda \cup \tilde\Lambda)$ and the partition $\Lambda \cup \tilde\Lambda=\bigcup\limits_{k}\Lambda_k'$, where $\Lambda'_{2k-1}=\Lambda_k$
     and $\Lambda'_{2k}=\{\tilde\lambda_k\}$.
\end{Remark}

\subsection{Necessary conditions.}
First we formulate and prove some elementary necessary conditions for the existence of a partition
of $\Lambda$ for which  $\mathcal{K}_p(\Lambda)$ is a
$\Big(\bigoplus\limits_{k=1}^{\infty}\ell^p_{n_k}\Big)_{\ell^1}$-frame.

\begin{Proposition}
    Let $\mathcal{K}_p(\Lambda)$ be a frame associated with $\Big(\bigoplus\limits_{k=1}^{\infty}\ell^p_{n_k}\Big)_{\ell^1}$. 
    Then the sequence $\{n_k\}_k$ is unbounded.
\end{Proposition}

\begin{proof}
 Suppose that 
 $n_k \leq N$ for all $k$ and put $A_l=\{z\in \D:1-2^{-l}\le|z| < 1-2^{-l-1}\}$.
 Let $n, L\in \mathbb{N}$ and let us estimate the frame sum for $f(z)=z^n$:
  $$\sum\limits_{\Lambda_k}(1-|\lambda|)|\lambda|^{qn}=
    \sum\limits_{l=0}^{\infty}\sum\limits_{\Lambda_k \cap A_l}
    (1-|\lambda|)|\lambda|^{qn}=
    \bigg(\sum\limits_{l=0}^{L}\sum\limits_{\Lambda_k \cap A_l}+
   \sum\limits_{l=L+1}^{\infty}\sum\limits_{\Lambda_k \cap A_l}\bigg)
   (1-|\lambda|)|\lambda|^{qn}.$$
    The second sum is bounded by $ N\sum\limits_{l=L+1}^{\infty}2^{-l}$. The first sum does not exceed
     $$(1-2^{-L-1})^{qn}
    {\sum\limits_{l=0}^{L}\sum\limits_{\Lambda_k \cap A_l}(1-|\lambda|)}\le (1-2^{-L-1})^{qn}
    \sum\limits_{\Lambda_k}(1-|\lambda|) \leq
    B^q (1-2^{-L-1})^{qn}.$$

   Choose $L$ such that the second sum is less than $A^q/4$ and then $n$ such that
   the first sum is less than $A^q/4$, which contradicts the lower frame
   inequality.
\end{proof}

\textbf{Proof of necessity in Theorem 1.4.}
\begin{proof}
We argue by contradiction. Suppose that for every $C,c>0$, there exist an
angle $\Gamma$ and $r\in(0,1)$ such that every
$\Lambda'\subset\Lambda\cap\Gamma\cap A_{r,1}$ satisfies either
\begin{equation}
 \sum_{\lambda\in\Lambda'}(1-|\lambda|) < c|\Gamma|
 \quad\text{or}\quad \text{Carl}(\Lambda') > C.
 \label{eq:failure}
\end{equation}
Since the partition $\Lambda= \bigcup\limits_{k=1}^{\infty}\Lambda_k$ generates a frame, there exist $A,B>0$ such that the following inequalities hold for every $f \in H^q$:
\[
 A^q\|f\|_q^q\le
 \sup\limits_{k}\sum_{\lambda\in\Lambda_k}(1-|\lambda|)|f(\lambda)|^q
 \le B^q\|f\|_q^q.
\]
In particular, $\text{Carl}(\Lambda_k)\le B_1$ for some constant $B_1$ and every $k$.
Apply \eqref{eq:failure} with the constants $B_1,\varepsilon$. For each $\Lambda_k$,
either $\text{Carl}(\Lambda_k \cap\Gamma\cap A_{r,1})>B_1$, which is not the case, or
\begin{equation}
 \sum_{\lambda\in\Lambda_k\cap\Gamma\cap A_{r,1}}(1-|\lambda|)
 < \varepsilon |\Gamma|,\quad k\in\N.
 \label{eq:mass}
\end{equation}

Next, let
\[
 \Gamma_1= \Gamma \cap \T=\{e^{it}:t\in(t_0-\tau,t_0+\tau)\},\qquad
 \Gamma_2=\left\{e^{it}:t\in
 \left(t_0-\frac{\tau}{2},t_0+\frac{\tau}{2}\right)\right\}.
\]
Thus $\Gamma_2$ is the concentric arc of half the length,
$|\Gamma_2|=|\Gamma_1|/2$.
Define $u:\T\to(0,\infty)$ by
\[
 u=
 \begin{cases}
  \displaystyle \left(\frac{4\pi}{|\Gamma_1|}\right)^{1/q}, &\text{on }\Gamma_2,\\[6pt]
  e^{-\Delta}, &\text{on }\T\setminus\Gamma_1,\\[3pt]
 \end{cases}
\]
where $\Delta>0$ is a parameter to be chosen later, $u \in C^1(\T)$, and
 $e^{-\Delta}\le u\le \displaystyle \left(\frac{4\pi}{|\Gamma_1|}\right)^{1/q}$.
Let $F$ be an outer function:
\[
 F(z)=\exp\left\{\int_{\T}\frac{z+\zeta}{\zeta-z}
 \log u(\zeta)\,dm(\zeta)\right\},\qquad z\in\D.
\]
Then $|F|=u$ almost everywhere on $\T$, and $|F|\in C(\overline\D)$. Moreover,
\[
 \|F\|_q^q=\int_\T u^q\,dm
 \ge\int_{\Gamma_2}u^q\,dm
 \ge\frac{4\pi}{|\Gamma_1|}\cdot\frac{|\Gamma_2|}{2\pi}=1.
\]
Hence $\|F\|_q\ge1$. The maximum modulus principle also gives
$|F(z)|\le\left(\frac{4\pi}{|\Gamma_1|}\right)^{1/q}$ in $\D$.
Set $f=z^{\nu/q}F$, where $\nu/q$ is a positive integer parameter to be fixed at the
end. Let us estimate the sum over a single block $\Lambda_k$.
We split the sum as follows:
\begin{align*}
 \sum_{\lambda\in\Lambda_k}(1-|\lambda|)|f(\lambda)|^q
 =\sum_{\Lambda_k\cap\Gamma\cap A_{r,1}}
       (1-|\lambda|)|f(\lambda)|^q+\sum_{\Lambda_k\cap(\Gamma\cap A_{r,1})^c}
       (1-|\lambda|)|f(\lambda)|^q.
\end{align*}
For the first term, $|f|\le|F|\le \left(\frac{4\pi}{|\Gamma_1|}\right)^{1/q}$ throughout
$\Gamma$, so \eqref{eq:mass} implies
\[
 \sum_{\Lambda_k\cap\Gamma\cap A_{r,1}}(1-|\lambda|)|f(\lambda)|^q
 \le\frac{4\pi}{|\Gamma|}
     \sum_{\Lambda_k\cap\Gamma\cap A_{r,1}}(1-|\lambda|)
 \le\frac{4\pi}{|\Gamma|}\, \varepsilon|\Gamma|=4\pi\varepsilon.
\]
On the other hand, $\Lambda_k\cap(\Gamma\cap A_{r,1})^c \subset (\Lambda_k \cap \{|\lambda|\le r'\}) \cup (\Lambda_k \cap A_{r',1}\cap \Gamma^c)$, where $r'\in(r,1)$
is a number such that $|F(z)|\le2e^{-\Delta}$ for
$z\in A_{r',1}\cap\Gamma^c$. Such an $r'<1$ exists by continuity of $|F|$
and the equality $|F|=e^{-\Delta}$ on $\T\setminus\Gamma$. Then
\[
 \sum_{\Lambda_k\cap \{|\lambda| \le r'\}}(1-|\lambda|)|f(\lambda)|^q
 \le (r')^\nu \sum_{\Lambda_k}(1-|\lambda|)|F(\lambda)|^q
 \le (r')^\nu B^q\|F\|_q^q,
\]
and
\[
 \sum_{\Lambda_k\cap A_{r',1}\cap\Gamma^c}(1-|\lambda|)|f(\lambda)|^q
 \le2^qe^{-q\Delta}\sum_{\Lambda_k\cap A_{r',1}\cap\Gamma^c}(1-|\lambda|)
 \le 2^q e^{-q\Delta}B^q.
\]
 Since $\|F\|_q^q \ge 1$, we have
\begin{align*}
 \sum_{\lambda\in\Lambda_k}(1-|\lambda|)|f(\lambda)|^q
 &\le 4\pi\varepsilon+2^qe^{-q\Delta}B^q
     +(r')^\nu B^q\|F\|_q^q\\
 &\le\bigl(4\pi\varepsilon+2^qe^{-q\Delta}B^q+(r')^\nu B^q\bigr)\|F\|_q^q.
\end{align*}

Note that $\|F\|_q = \|f\|_q$. Choosing sufficiently small $\varepsilon$ and sufficiently large $\Delta$
and $\nu$, we obtain
 a contradiction with the lower frame inequality.
\end{proof}

As an immediate corollary of this part of the theorem, we obtain the following necessary conditions. 

\begin{Corollary}
    Let the system $\mathcal{K}_p(\Lambda)$ be a frame associated with $\Big(\bigoplus\limits_{k=1}^{\infty}\ell^p_{n_k}\Big)_{\ell^1}$. 
    Then for every $r\in (0,1)$ and every angle $\Gamma$ the following condition
    holds: $$\sum\limits_{\Lambda \cap A_{r,1} \cap \Gamma}(1-|\lambda|)=\infty.$$
\end{Corollary}

\noindent \textbf{Problem.} It follows from Corollary 3.4 that if the system $\mathcal{K}_p(\Lambda)$ is
an $X_p$-frame, then $\text{Clos }\Lambda \supset \T$. 
We have already seen that this may not hold for an arbitrary representing system; see Corollary 2.4. Assume now that a system $\mathcal{K}_p(\Lambda)$ is representing for $H^p$ and
$\Lambda$ has no accumulation points in $\mathbb{D}$. An intriguing question is whether in this case we still have $\text{Clos }\Lambda \supset \T$.

\subsection{Sufficient conditions.}

\textbf{Proof of sufficiency in Theorem 1.4.}
\begin{proof}
Let $\{{s_k\}}_{k=1}^{\infty}$ be a sequence of positive integers tending to infinity and let $r_1 \in (1-2\pi/s_1,1)\cap(0,1)$. We construct $\Lambda_k$ and radii $r_k\to1-$ inductively.
Set
\[
 \Gamma_{k,j}=\left\{re^{it}:r\in(0,1),\quad
 t\in\left[\frac{2\pi j}{s_k},\frac{2\pi(j+1)}{s_k}\right)\right\},
 \qquad j=0,\ldots,s_k-1.
\]
For each fixed $k$ and each sector $\Gamma_{k,j}$, the hypothesis yields a set
 $\Lambda'_{k,j}$ satisfying
\[
 \Lambda'_{k,j}\subset\Lambda\cap\Gamma_{k,j}\cap A_{r_k,1},\qquad
 \text{Carl}(\Lambda'_{k,j})\le C,\qquad
 \sum_{\lambda\in\Lambda'_{k,j}}(1-|\lambda|)\ge c|\Gamma_{k,j}|,
\]
where $|\Gamma_{k,j}|=2\pi/s_k$. Then there exist a radius $r_{k,j}$ and a finite set $\Lambda_{k,j}'' \subset \Lambda_{k,j}'$ such that
$$
 \sum\limits_{\Lambda''_{k,j}\cap A_{r_k,r_{k,j}}}(1-|\lambda|)
 \ge c_1|\Gamma_{k,j}|,$$
 
  where $c_1=c/2$. Then we put
\[
 \Lambda_{k,j}:=\Lambda''_{k,j}\cap A_{r_k,\,r_{k,j}},\qquad
 r_{k+1}:=\max_j\left\{r_{k,j},\ 1-\frac{2\pi}{s_{k+1}}\right\}, \qquad
 \Lambda_k:=\bigcup_{j=0}^{s_k-1}\Lambda_{k,j},
\]
so that $\text{Carl}(\Lambda_{k,j})\le C$.

\textbf{Upper frame inequality}. Let
\[
 \mu_k=\sum_{\lambda \in \Lambda_k}(1-|\lambda|)\delta_\lambda,\qquad
 \mu_{k,j}=\sum_{\lambda \in\Lambda_{k,j}}(1-|\lambda|)\delta_\lambda.
\]
Now we estimate $\text{Carl}(\Lambda_k)$.
Let $\zeta\in\T$ and $h\le 2\pi/s_k$. Then $S(\zeta,h)$ intersects at most two sectors
$\Gamma_{k,j}$, say $\Gamma_{k,j_1}$ and $\Gamma_{k,j_2}$. Thus
\begin{align*}
 \mu_k(S(\zeta,h))
 &\le\mu_k\bigl(S(\zeta,h)\cap\Gamma_{k,j_1}\bigr)
     +\mu_k\bigl(S(\zeta,h)\cap\Gamma_{k,j_2}\bigr)\\
 &=\mu_{k,j_1}(S(\zeta,h))+\mu_{k,j_2}(S(\zeta,h))\le2Ch.
\end{align*}
Now take $j=1,\ldots,s_k-1$ and
$h\in[2\pi j/s_k,\,2\pi(j+1)/s_k)$.
Then $S(\zeta,h)$ intersects at least $j$ and at most $j+2$ sectors
$\Gamma_{k,l}$, say
$\Gamma_{k,l_1},\ldots,\Gamma_{k,l_{j+2}}$. Hence
\begin{align*}
 \mu_k(S(\zeta,h))
 \le\sum_{s=1}^{j+2}\mu_k\bigl(S(\zeta,h)\cap\Gamma_{k,l_s}\bigr)
  =\sum_{s=1}^{j+2}\mu_{k,l_s}(S(\zeta,h))\\
 \le \sum_{s=1}^{j+2}\mu_{k,l_s}
       \left(\Gamma_{k,l_s}\cap A_{1-\frac{2\pi}{s_k},\,1}\right) \le(j+2)\cdot\frac{2\pi}{s_k}\cdot C\le3Ch.
\end{align*}
Therefore $\text{Carl}(\Lambda_k)\le 3C$, which yields the upper frame inequality with some constant $C_1$.

\textbf{Lower frame inequality for $f\in A(\D)$.} Set
\[
 I_{k,j}=\left\{e^{it}: t\in\left[\frac{2\pi j}{s_k},\frac{2\pi(j+1)}{s_k}\right)\right\},
 \qquad j=0,\ldots,s_k-1.
\]By the Minkowski inequality

    $$\Big(\sum\limits_{\lambda \in\Lambda_k}(1-|\lambda|)|f(\lambda)|^q\Big)^{1/q}= \Big(\sum\limits_{j=0}^{s_k-1}\sum\limits_{\lambda \in \Lambda_{k,j}}(1-|\lambda|)|f(\lambda)|^q\Big)^{1/q}$$
    $$=\Bigg(\sum\limits_{j=0}^{s_k-1}\sum\limits_{\lambda \in \Lambda_{k,j}}(1-|\lambda|)|f(\lambda)|^q\int\limits_{I_{k,j}}\frac{2 \pi dm(\zeta)}{|I_{k,j}|}\Bigg)^{1/q}\ge 
    \Bigg(\sum\limits_{j=0}^{s_k-1}\sum\limits_{\lambda \in \Lambda_{k,j}}(1-|\lambda|)\int\limits_{I_{k,j}}|f(\zeta)|^q\frac{2 \pi dm(\zeta)}{|I_{k,j}|}\Bigg)^{1/q}$$
    $$-\Bigg(\sum\limits_{j=0}^{s_k-1}\sum\limits_{\lambda \in \Lambda_{k,j}}(1-|\lambda|)\int\limits_{I_{k,j}}|f(\zeta)-f(\lambda)|^q\frac{ 2\pi dm(\zeta)}{|I_{k,j}|}\Bigg)^{1/q}.$$

    For every $\varepsilon >0$ and sufficiently large $k$ (depending on $\varepsilon$), the last sum is bounded by 
    $\varepsilon \Big(\sum\limits_{j=0}^{s_k-1}\sum\limits_{\lambda \in \Lambda_{k,j}}(1-|\lambda|)\Big)^{1/q}\le \varepsilon C_1$ by uniform continuity of $f$ in $\overline{\D}$.

   For the first sum we get 
   $$\Bigg(\sum\limits_{j=0}^{s_k-1}\sum\limits_{\lambda \in \Lambda_{k,j}}(1-|\lambda|)\int\limits_{I_{k,j}}|f(\zeta)|^q\frac{2 \pi dm(\zeta)}{|I_{k,j}|}\Bigg)^{1/q}=
   \Bigg(\sum\limits_{j=0}^{s_k-1}\int\limits_{I_{k,j}}|f(\zeta)|^q\frac{2 \pi dm(\zeta)}{|I_{k,j}|}\sum\limits_{\lambda \in \Lambda_{k,j}}(1-|\lambda|)\Bigg)^{1/q}$$
   $$\ge \bigg( 2\pi c_1\sum\limits_{j=0}^{s_k-1}\int\limits_{I_{k,j}}|f(\zeta)|^qdm(\zeta) \bigg)^{1/q}=(2\pi c_1)^{1/q}\|f\|_q.$$

   Taking $\varepsilon $  small enough, we obtain the lower estimate for sufficiently large $k$:
   $$\Big(\sum\limits_{\lambda \in \Lambda_k}(1-|\lambda|)|f(\lambda)|^q\Big)^{1/q}\ge c_2\|f\|_q$$
   with some constant $c_2>0$.

   \textbf{Lower frame inequality for an arbitrary function.} Now we obtain the lower frame inequality for arbitrary $f \in H^q$. For such an $f$, choose $g \in A(\D)$ such that $\|g-f\|_{q}\ \le \varepsilon$.
    Then by the Minkowski inequality and the upper frame inequality we get 
    \begin{align*}\bigg(\sum\limits_{\Lambda_k}(1-|\lambda|)|f(\lambda)|^q\bigg)^{1/q} &\ge
     \bigg(\sum\limits_{\Lambda_k}(1-|\lambda|)|g(\lambda)|^q\bigg)^{1/q}-
     \bigg(\sum\limits_{\Lambda_k}(1-|\lambda|)|f(\lambda)-g(\lambda)|^q\bigg)^{1/q} \\
      &\ge \bigg(\sum\limits_{\Lambda_k}(1-|\lambda|)|g(\lambda)|^q\bigg)^{1/q}-C_1\|f-g\|_{q}.
      \end{align*}

   Consequently,
   \begin{align*}
     \sup\limits_{k}\bigg(\sum\limits_{\Lambda_k}(1-|\lambda|)|f(\lambda)|^q\bigg)^{1/q}\ge
    \sup\limits_{k}\bigg(\sum\limits_{\Lambda_k}(1-|\lambda|)|g(\lambda)|^q\bigg)^{1/q} -
     \varepsilon C_1\\
   \ge c_2\|g\|_q-\varepsilon C_1 \ge
    c_2(\|f\|_q-\|f-g\|_q)-\varepsilon C_1 \ge c_2\|f\|_q-C_2\varepsilon.
   \end{align*} 
   Taking $\varepsilon$ small enough, we get the lower frame inequality for $f$.

   Thus, the blocks $\{\Lambda_{k}\}_k$ generate a frame. By Remark 3.2 the whole 
   set $\Lambda$ admits a partition with the required frame property.   
\end{proof}

As an immediate corollary of this part of the theorem, we obtain the following sufficient condition.

\begin{Corollary}
    Suppose that for some constants $c, C>0$ there exist sequences $r_k \to 1-$, $r'_k\in (r_k,1)$ such that for any angle $\Gamma$ with $|\Gamma|\ge (1-r'_k)$ the following
        condition holds:
        $$c|\Gamma| \le\sum\limits_{\Lambda\cap A_{r_k,r'_k}\cap \Gamma}(1-|\lambda|)\leq C|\Gamma|.$$ 
        Then the system $\mathcal{K}_p(\Lambda)$ is a $\Big(\bigoplus\limits_{k=1}^{\infty}\ell^p_{n_k}\Big)_{\ell^1}$-frame with respect to some partition.
\end{Corollary}

The frames constructed in \cite{ST18} and \cite{BB23} have the structure described
in the preceding corollary.

Now we present an example of an $X_p$-frame with a different structure. 

\begin{Proposition}
     Let $p \in(1,\infty)$ and let $\Lambda=\{\lambda_n=(1-1/n)e^{i \log n}: n \ge 1\}$. Then there exists a partition $\Lambda=\bigcup\limits_{k}\Lambda_{k}$
      such that $n_k=\text{card }\Lambda_k < \infty$ and the corresponding system
  $\mathcal{K}_p(\Lambda)$ is a frame for $H^p$ with respect to the space
  $\Big(\bigoplus\limits_{k}\ell^p_{n_k} \Big)_{\ell^1}$.
\end{Proposition}
\begin{proof}
    Let $r \in (0,1)$, $0 < \beta-\alpha \le 2\pi$, $\Gamma=\{z\ \in \C: \text{arg }z \in (\alpha, \beta)\}$. Then $\lambda_n \in \Gamma \cap A_{r,1}$ if
    $n \ge (1-r)^{-1}$ and $n \in (e^{\alpha+2\pi k}, e^{\beta+2\pi k})$ with some positive integer $k$ such that $e^{\alpha+2\pi k}> (1-r)^{-1}$. Let $N_k$ denote the set of all such integers $n$, and put
    $\Lambda'_k=\{\lambda_n\}_{n \in N_k}$. Then $$ \sum\limits_{\Lambda_k'}(1-|\lambda|)= \sum\limits_{n=[e^{\alpha+2\pi k}]+1}^{[e^{\beta+2\pi k}]}\frac{1}{n}=\beta-\alpha+o(1).$$

    Next, let $\zeta=e^{i\tau}$, $h \in (0,2 \pi]$. If $h<e^{-\beta-2\pi k}$, then $S(\zeta,h)\cap \Lambda_k' =\varnothing$. Moreover,  
    we have $$S(\zeta,h)\cap \Lambda'_k \subset \{\lambda_n: e^{\tau-h/2+2\pi k}\le n \le e^{\tau+h/2+2\pi k}\}.$$
    Thus, we obtain 
    \begin{align*}
    \mu_{\Lambda_{k}'}S(\zeta,h)&\le \sum\limits_{n=[e^{\tau-h/2+2\pi k}]+1}^{[e^{\tau+h/2+2\pi k}]} \frac{1}{n} \le
    \int\limits_{[e^{\tau-h/2+2\pi k}]}^{e^{\tau+h/2+2\pi k}}\frac{dx}{x} \\
        &\le \tau+h/2+2\pi k- \log(e^{\tau-h/2+2\pi k}-1) \\
        &= h-\log(1-e^{h/2-\tau-2\pi k}).
    \end{align*}
    For $h\ge e^{-\beta-2\pi k}$, we have 
    $$-\log(1-e^{h/2-\tau-2\pi k}) \le \frac{e^{h/2-\tau-2\pi k}}{1-e^{h/2-\tau-2\pi k}} \le Ch$$ with some absolute constant $C$. The conclusion follows from Theorem 1.4.
\end{proof}

\section{Representing systems of dilations of a function for $H^2$}

Let $g(z)=\sum\limits_{n=0}^{\infty}\beta_n z^n$ be an analytic function in $\D$, $\beta_n \neq 0$,
    $\lim\limits_{n \to \infty} |\beta_n|^{1/n}=1$. We will present a sufficient condition on
    $\Lambda \subset \D$ for the system
    $\{g(\bar\lambda z)\}_{\lambda \in \Lambda}$ to be a representing system for $H^2$.
    The following lemma, which follows by a direct computation, gives an integral representation of a function
    $f \in H^2$ in terms of $g$.
    
\begin{Lemma}
    If $f(z)=\sum\limits_{n=0}^{\infty} a_nz^n \in H^2$, 
    $F_{\rho}(z)=\sum\limits_{n=0}^{\infty}\rho^n a_n z^n/ \beta_n$, and $0<\rho<1$, then 
    $$f(rz)=\int\limits_{\T}F_{r/R}(\zeta) {g(R\bar\zeta z)}dm(\zeta), \quad 0<r<R<1.$$
\end{Lemma} 

\begin{Remark}
    Notice that $$\bigg(\displaystyle \int_\T  |F_{\rho} (\zeta)| dm(\zeta)\bigg)^2 \le
    \int_\T  |F_{\rho} (\zeta)|^2 dm(\zeta)= \sum_{n=0}^\infty |a_n|^2 \rho^{2n}/|\beta_n|^2 
    \le \omega_1(\rho)\|f\|_2^2,$$
    where $\omega_1(r)= \sup\limits_{n\ge 0} (r^{2n}/|\beta_n|^2)$.
\end{Remark} 

We introduce two more characteristics of the function $g$:
$$\omega_2(r)= \sum\limits_{n=1}^{\infty} n^2 |\beta_n|^2 r^{2n},
\qquad \omega_3(r)=\sum\limits_{n=0}^{\infty} |\beta_n|^2r^{2n}.$$

Let $k\in \N$, $R_k \in(0,1)$, $N_k, M_k \in \N$
and $R_{k,l} \in (0,1)$, $l=1, \dots, M_k$ with $R_{k,1} <R_{k,2} \dots 
< R_{k, M_k} =  R_k$. Put
$$
\Lambda = \{w_{k, l, j} :\ 1\le l \le M_k, 1\le j\le N_k, k\in\N\},  \qquad  w_{k, l, j} = R_{k,l} \exp \Big(\frac{2\pi i j}{N_k}\Big).
$$

\begin{Theorem}
  Let $0< \varepsilon<1/(4\pi^2)$, $R_{k, 1} \to 1$, and
$$
\omega_1(R_k) \omega_2(R_k)/N_k^2 \le \varepsilon, \qquad
\omega_1(R_k) \omega_3(R_k) = O(M_k^2)
$$
as $k\to\infty$. Then $\{g(\bar\lambda z)\}_{\lambda\in\Lambda}$ is a representing system for $H^2$.
\end{Theorem} 

\begin{proof}
    Using Lemma 4.1 with $r_k=R_{k,1}^3<R_{k,1}$, we obtain
$$
f_{r_k}(z)=f(r_k z) = \frac{1}{M_k} \sum_{l=1}^{M_k} \int_ \T F_{r_k/R_{k,l}} (\zeta) g(R_{k,l}\bar\zeta z) dm(\zeta).
$$
Consider the following discretization:
$$
S_k(z) = \frac{1}{M_k} \sum_{l=1}^{M_k} \sum_{j=1}^{N_k} 
\bigg(\int_{I_{k,j}} F_{r_k/R_{k,l}} (\zeta) dm(\zeta)\bigg) g(\overline{w}_{k,l,j} z),
$$
where $I_{k,j} = \left\{e^{it}:  t\in \left [\frac{(2 j -1)\pi }{N_k},\frac{(2 j +1)\pi }{N_k}\right]\right\}$. 

Now we describe the first step of approximation. For the moment we fix the number $k$ and omit
 it from the notation,
 so we write $r, R, R_l, w_{l,j}, M, N$ instead of 
$r_k, R_k, R_{k,l}, w_{k,l,j}, M_k, N_k$. 
Recall that $R_1<\dots < R_M = R$ and suppose $r= R^3_1 < R_1^2$.

We get
\begin{align}
f(rz) - S(z)  =  \frac{1}{M} \sum_{l=1}^{M} \sum_{j=1}^{N} 
\int_{I_j}  F_{r/R_{l}} (\zeta) \big (g( R_{l}\bar\zeta z) - g(\overline{w}_{l,j}z)\big)
dm(\zeta) = \sum_{n=1}^\infty c_n z^n,
\label{eq:1}
\end{align}
where
\begin{align}
    c_n =  \frac{\beta_n}{M} \sum_{l=1}^{M} \sum_{j=1}^{N} 
\int_{I_j}  \big( (R_l\bar \zeta)^n - \overline{w}_{l,j}^n \big) F_{r/R_{l}}(\zeta) dm(\zeta).
\end{align} 
Since $w_{l,j} = R_l \zeta_j$, $\zeta_j = e^{2\pi i j/N}  \in I_j$, for every $ \zeta \in I_j$ we get
\begin{align}
|(R_l \zeta)^n - w_{l,j}^n| =R_l^n \cdot |\zeta - \zeta_j| \cdot  \Big| \sum_{s=0}^{n-1} \zeta^s \zeta_j^{n-1-s}\Big|
\le \frac{2\pi n R_l^n}{N} \le \frac{2 \pi n R^n}{N}.
\end{align} 
Therefore,
\begin{align}
|c_n| \le \frac{2\pi n |\beta_n|  R^n }{M N} \sum_{l=1}^{M}  \int_\T |F_{r/R_{l}}(\zeta)| dm(\zeta).
\end{align} 
Recall that $r\le R_1^2$, which implies $r/R_l \le R$.
Consequently, 
\begin{align}
\int_\T |F_{r/R_{l}}(\zeta)| dm(\zeta) \le (\omega_1(R))^{1/2} \|f\|_2.
\end{align}
We conclude that
\begin{align}
|c_n| \le \frac{2\pi n |\beta_n| R^n (\omega_1(R))^{1/2}}{N} \|f\|_2.
\label{eq:2}
\end{align}
Therefore, 
$$
\|f_r - S\|_2^2  = \sum_{n=1}^\infty |c_n|^2 \le 
\frac{4\pi^2\omega_1(R)   \|f\|^2_2}{N^2} \sum_{n=1}^\infty n^2 |\beta_n|^2 R^{2n} = 
\frac{4\pi^2\omega_1(R)  \omega_2(R)}{N^2} \|f\|^2_2.
$$
Thus, $\|f_r-S\|_2 \le 2\pi \sqrt{\varepsilon} \|f\|_2$. 

Now we show that the norms of all intermediate sums are uniformly bounded by $C\|f\|_2$. Put 
$$
T_{\tilde M}(z) = \frac{1}{M} \sum_{l=1}^{\tilde M} \sum_{j=1}^{N} 
\bigg( \int_{I_j}  F_{r/R_{l}} (\zeta) dm(\zeta) \bigg) g(\overline{w}_{l,j}z) 
$$
and 
$$
T_{\tilde M, \tilde N}(z) = \frac{1}{M} \sum_{j=1}^{\tilde N} 
\bigg( \int_{I_j}  F_{r/R_{\tilde M}} (\zeta) dm(\zeta) \bigg) g(\overline{w}_{\tilde{M},j}z),
$$
where $1\le \tilde M\le M$ and $1\le \tilde N \le N$. 

First, we consider the sum $T_{\tilde M}(z)$ over complete circles.
Repeating the previous estimates \eqref{eq:1}--\eqref{eq:2} and using Lemma 4.1, we conclude
that for every $1\le \tilde M \le M$
$$
\bigg\| T_{\tilde M}(z)  -
\frac{1}{M} \sum_{l=1}^{\tilde M} \int_ \T  F_{r/R_{l}} (\zeta) g(R_l \bar\zeta z) dm(\zeta) 
\bigg\|_2^2
= \bigg\| T_{\tilde M}(z) - \frac{\tilde M}{M} f(rz) \bigg\|_2^2$$
$$\le  \frac{1}{M}\sum\limits_{l=1}^{\tilde M}\bigg \| \sum_{j=1}^{N} 
\bigg( \int_{I_j}  F_{r/R_{l}} (\zeta) dm(\zeta) \bigg) g(\overline{w}_{l,j}z) - f(rz) \bigg \|^2_2$$
$$ \le \frac{\tilde M}{M}  \frac{4\pi^2 \omega_1(R)  \omega_2(R)}{N^2} \|f\|^2_2
\le \frac{4\pi^2 \omega_1(R)  \omega_2(R)}{N^2} \|f\|^2_2 .$$
Therefore, $\| T_{\tilde M} \|_2 \le C\|f\|_2$ with some constant $C$.

Now we have to estimate the sum $T_{\tilde M, \tilde N}(z)$ over an incomplete circle.
Again, by the previous estimates of $f_r-S$ with $I = \cup_{j=1}^{\tilde N} I_j$ instead of $\T$,
$$
\bigg\|  T_{\tilde M, \tilde N}(z) -  
\frac{1}{M} \int_ I F_{r/R_{\tilde M}} (\zeta) g(R_{\tilde M} \bar\zeta z) dm(\zeta) 
 \bigg\|_2^2 \le  \frac{4\pi^2\omega_1(R)  \omega_2(R)}{N^2} \|f\|^2_2.
$$
Using the Taylor series of $g$, we obtain  
$$
\int_ I  F_{r/R_{\tilde M}} (\zeta) g(z R_{\tilde M} \bar\zeta) dm(\zeta) 
= \sum_{n=0}^\infty  d_n R_{\tilde M}^n \beta_n z^n, 
$$
where
$$
d_n = \int_ I  F_{r/R_{\tilde M}} (\zeta) \bar \zeta^n dm(\zeta).
$$

Estimating $d_n$ (see Remark 4.2) and using the inequality $r/R_{\tilde M} \le R$, we conclude that
$$
\bigg\| \frac{1}{M}\int_ I  F_{r/R_{\tilde M}} (\zeta) g(z R_{\tilde M} \bar\zeta) dm(\zeta) 
\bigg\|_2^2  \le \frac{\omega_1(R) \|f\|^2_2}{M^2} \sum_{n=0}^\infty R^{2n}|\beta_n|^2$$ 
$$= \frac{ \omega_1(R)\omega_3(R)}{M^2}  \|f\|^2_2 \le   C \|f\|^2_2$$
with some constant $C$.

Now we can finish the proof. We fix $\sigma >0$ such that $\gamma =\sigma +2\pi \sqrt{\varepsilon}<1$ and
for arbitrary $f\in H^2$ choose $r_{k_1}$ such that $\|f_{r_{k_1}}-f\| \le\sigma \|f\|$. 
 Then we get the function  
$$
f_1(z) = f(z) -  S_{k_1}(z)
$$
for which the estimate $\|f_1\| \le \gamma \|f\|$ holds.

Proceeding in this way, we find a sequence of indices $k_p$ and coefficients $c_{k_p, l, j}$ such that 
$$
\|f_J\| = \bigg\|f- \sum\limits_{p=1}^{J}\sum_{l=1}^{M_{k_{p}}} \sum_{j=1}^{N_{k_p}}
c_{k_p,l,j}g(\overline{w}_{k_p,l,j}z)\bigg\| \le \gamma^J \|f\|. 
$$

Now we only need to prove that the series 
$$
\sum\limits_{p=1}^{\infty}\sum_{l=1}^{M_{k_{p}}} \sum_{j=1}^{N_{k_p}}
c_{k_p,l,j}g(\overline{w}_{k_p,l,j}z)
$$
converges to $f$ in $H^2$. Indeed, for every $J\in \mathbb{N}$ 
and $1\le \tilde M\le M_{k_J}$, $1\le n\le N_{k_J}$
we get
$$
\begin{aligned}
&\bigg\|f - \sum_{p=1}^{J-1} \sum_{l=1}^{M_{k_p}}
 \sum\limits_{j=1}^{N_{k_p}}c_{k_p,l,j}g(\overline{w}_{k_p,l,j}z)  - \sum_{l=1}^{\tilde M-1}
 \sum\limits_{j=1}^{N_{k_J}}c_{k_J,l,j}g(\overline{w}_{k_J,l,j}z) -
\sum\limits_{j=1}^{n}c_{k_J,\tilde M,j}g(\overline{w}_{k_J,\tilde M,j}z)\bigg\| \\
& \le \bigg\|f - \sum_{p=1}^{J-1} \sum_{l=1}^{M_{k_p}}
 \sum\limits_{j=1}^{N_{k_p}}c_{k_p,l,j}g(\overline{w}_{k_p,l,j}z) \bigg\| +
\bigg\| \sum_{l=1}^{\tilde M-1}
 \sum\limits_{j=1}^{N_{k_J}}c_{k_J,l,j}g(\overline{w}_{k_J,l,j}z)\bigg\| \\
 &+\bigg\| \sum\limits_{j=1}^{n}c_{k_J,\tilde M,j}g(\overline{w}_{k_J,\tilde M,j}z)\bigg\| 
 \le  (1+C_1+C_2) \|f_{J-1}\| \le (1+C_1+C_2)  \gamma^{J-1} \|f\| \to 0
\end{aligned}
$$
as $J\to\infty$. This completes the proof.
\end{proof}

\begin{Remark}
Using the technique presented in this section, one can construct representing systems of dilations of a single function in $H^p$, $p \in(1,\infty)$.
\end{Remark}
\section{A linear method for the choice of coefficients}
\subsection{A linear method for $H^2$}
\begin{definition}
    We say that a generalized frame $(F, \mathcal{F}, X)$ \textit{allows a linear method for
     the choice of coefficients} if there exists 
    a sequence of functionals $\{g_n\}_n \subset F^*$ such that for every $f \in F$ the following two conditions hold:
    $$(1) \quad \{g_n(f)\}_n \in X , \qquad (2) \quad f=\sum\limits_{n}g_n(f)f_n.$$
\end{definition}

\begin{Remark}
    Let a system $\mathcal{F}=\{f_n\}_{n=1}^{\infty} \subset F\setminus\{0\}$
     be a representing system for a Banach space $F$, and let $X$ be its coefficient space:  
    $$X=\Bigg\{x=\{x_n\}_{n=1}^{\infty}: \sum\limits_{n \ge 1} x_n f_n \text{ converges },
    \, \|x\|_X= \sup\limits_{N\ge 1} \bigg\| \sum\limits_{n=1}^{N} x_n f_n \bigg\| \Bigg \}.$$
    Then the triple $(F, \mathcal{F}, X)$ is a generalized frame; see \cite{Ter10}.    
\end{Remark}

The following proposition was proved by P. A. Terekhin in \cite{Ter10}.

\begin{Proposition} Let $(F, \mathcal{F}, X)$ be a generalized frame. Then the following statements are equivalent:

(1) the frame $(F, \mathcal{F}, X)$ allows a linear method; 

(2) there exists a Banach space $F_1$, a bounded linear isomorphism $J:X \to F_1$ and a continuous projection $P:F_1 \to F$
such that $S=PJ$, where $S:X \to F$ is the synthesis operator 

$$Sx=\sum\limits_{n}x_n f_n.$$
\end{Proposition}

\begin{Corollary}
 If a triple 
 $(H^2, \mathcal{K}_2(\Lambda), \Big(\bigoplus\limits_{k=1}^{\infty}\ell^2_{n_k}\Big)_{\ell^1})$ is a frame, then it does not allow a linear method.
\end{Corollary}

\begin{proof}
    The space $X_2=\Big(\bigoplus\limits_{k=1}^{\infty}\ell^2_{n_k}\Big)_{\ell^1}$ has the Schur property as an $\ell^1$-sum
    of finite-dimensional spaces (see \cite{Tan98}).
    If $(H^2, \mathcal{K}_2(\Lambda), X_2)$
    allows a linear method, then
    it follows from Proposition 5.3 that the space $F_1$
     has the Schur property as a space isomorphic
     to $X_2$, and $H^2$ has the Schur property as a subspace of $F_1$,
     which is not true, since $H^2$ is an infinite-dimensional Hilbert 
     space.   
\end{proof}

Now we give the proof of Theorem 1.6.

\begin{proof}
    Without loss of generality, we may assume that $n_k(1-r_k)\le n_{k+1}(1-r_{k+1})$  and
     $\sum\limits_{k=1}^{\infty} \frac{1}{n_{k}(1-r_{k})} <\infty$.
 
    Put $\zeta_{k,j}=\exp(2\pi i j / n_k)$, $I_{k,j}=\{e^{it}:t \in [2\pi  (j-1)/n_k, 2\pi  j/n_k]\}$ and
    also put $f_0=f$, $f_{k+1}(z)=f_k(z)-S_{k+1}(f_k,z)$, where $$S_k(f,z)=\sum\limits_{j=1}^{n_k} \int\limits_{I_{k,j}}f(\zeta)dm(\zeta) \cdot \frac{1}{1-r_k \bar{\zeta}_{k,j}z}.$$
    Put $T_n f(z)=f(z)-f(r_{n}z)$, $V_n f = f(r_n z)-S_n(f,z)$. Then $$f_{n+1}=T_{n+1}f_n+V_{n+1}f_n=(T_{n+1}+V_{n+1})...(T_1+V_1)f.$$ 
    Put $B_n=T_n+V_n$, $A_n=B_n...B_1$. Then $\|A_n\| \leq \prod\limits_{q=1}^{n}(1+\|V_q\|)$ as $\|T_k\|=1$. 
    
    \textbf{Estimate of the norm of $V_q$.}
    For the operator $V_q$ we obtain
     $$\displaystyle V_q f = \int\limits_{\T}\frac{f(\zeta)dm(\zeta)}{1-r_q \bar{\zeta}z}-\sum\limits_{j=1}^{n_q} 
    \int\limits_{I_{q,j}}\frac{f(\zeta)dm(\zeta)}{1-r_q \bar{\zeta}_{q,j}z}=\sum\limits_{j=1}^{n_q} \int\limits_{I_{q,j}} f(\zeta)\left(\frac{1}{1-r_q \bar{\zeta}z}-
    \frac{1}{1-r_q \bar{\zeta}_{q,j}z}\right) dm(\zeta)$$

    $$=\sum\limits_{j=1}^{n_q} \int\limits_{I_{q,j}} f(\zeta)\frac{ r_q(\bar\zeta- \bar\zeta_{q,j})z}
    {(1-r_q  \bar{\zeta}z)(1-r_q \bar{\zeta}_{q,j}z)} dm(\zeta).$$
    By a direct calculation 
    $$\Bigg|\frac{1}
    {(1-r_q  \bar{\zeta}z)(1-r_q \bar{\zeta}_{q,j}z)} \Bigg|=\Bigg|\frac{1-r_q \bar{\zeta}_{q,j}z+r_q \bar{\zeta}_{q,j}z -r_q  \bar{\zeta}z}
    {(1-r_q  \bar{\zeta}z)^2(1-r_q \bar{\zeta}_{q,j}z)} \Bigg| $$ 
    $$\le \bigg(1+\frac{2 \pi}{n_q (1-r_q)}\bigg)\frac{1}{|1-r_q  \bar{\zeta}z|^2}.$$
    Therefore
    $$|V_q f| \le \sum\limits_{j=1}^{n_q} \int\limits_{I_{q,j}} |f(\zeta)|\frac{|\zeta-\zeta_{q,j}|}
    {|1-r_q z \bar{\zeta}||1-r_q \bar{\zeta}_{q,j}z|} dm(\zeta) $$ 
    $$\le \frac{2\pi}{n_q} \bigg(1+\frac{2\pi}{n_q(1-r_q)}\bigg)\sum\limits_{j=1}^{n_q}
    \int\limits_{I_{q,j}} \frac{|f(\zeta)|}{|1-r_q z \bar \zeta|^2}dm(\zeta)$$
    $$=\frac{2\pi}{n_q} (1+a_q) \int\limits_{\T} \frac{|f(\zeta)|}{|1-r_q z \bar \zeta|^2}dm(\zeta),
    \quad \text{where } a_q=\frac{2\pi}{n_q(1-r_q)}.$$
    By H\"older's inequality we get
    $$
    \|V_q f\|_{H^2}^2 \le 
    \frac{(2\pi)^2}{n_q^2} (1+a_q)^2\int_\T \bigg( \int_{\T} \frac{|f(\zeta)|^2}{|1 - r_q \bar \zeta z|^2} dm(\zeta)\bigg) 
    \bigg( \int_{\T} \frac{dm(\zeta)}{|1 - r_q \bar \zeta z|^2} \bigg) dm(z).
    $$
    As $\int_{\T} |1 - r_q u|^{-2} dm(u) =(1-r_q^2)^{-1}$, 
    we conclude that
    $$
    \|V_q f\|_{2} \le a_q(1+a_q)\|f\|_2. 
    $$
    Thus, the operators $A_n$ are uniformly bounded. Put $a =  \prod\limits_{q=1}^{\infty}(1+a_q(1+a_q))$.
    Then $\sup\limits_{n}\|A_n\|\leq a$.

    We have to show that $A_n f \to 0$. Note that $A_nf=T_nA_{n-1}f+V_nA_{n-1}f$. For the second term we obtain 
    $$\|V_nA_{n-1}f\|\leq a \cdot a_n(1+a_n)\|f\| \to 0.$$

    To see that $T_nA_{n-1}f \to 0$, we have to use some commutation properties of the operators $T_n$ and $V_m$. Direct calculation shows that $T_mT_n=T_nT_m$.
    For the pair of operators $T_m$ and $V_k$ we have:
    $$T_m V_k f(z)= f(r_k z)-S_k(f, z)-(f(r_kr_m z)-S_k(f, r_m z)),$$
    $$V_kT_mf(z)=f(r_k z)-f(r_k r_m z)-(S_k(f,z)-S_k(f(r_m \cdot),z)).$$
    Put $$R_{k,m}f(z)= (T_m V_k - V_k T_m) f(z)=S_k(f, r_m z)-S_k(f(r_m \cdot),z).$$
    Then $$R_{k,m}f(z)=S_k(f, r_m z)-f(r_k r_m z)+f(r_k r_m z)-S_k(f(r_m \cdot),z)$$ 
    $$=-(V_kf)(r_mz)+(V_kf(r_m \cdot))(z).$$ 
    Therefore, $\|R_{k,m}\|\leq 2\|V_k \| \leq 2(1+a_k)a_k.$
    We also need the following recurrence relation:
    $$A_{N+2n}=B_{N+2n-1} T_{N+2n}A_{N+2n-2}+R_{N+2n-1, N+2n} A_{N+2n-2}+V_{N+2n}A_{N+2n-1},$$
    where $$\|R_{N+2n-1, N+2n} A_{N+2n-2}+V_{N+2n}A_{N+2n-1}\| \leq 3a \cdot a_{N+2n-1}(1+a_{N+2n-1}).$$
    Consequently, $$A_{N+2n}=\bigg(\prod\limits_{j=0}^{2n-2}B_{N+2n-1-j}\bigg)T_{N+2n}A_N+\hat{R}_{N,n},$$
    where $\|\hat{R}_{N,n} \|\leq 3 a^2 \sum\limits_{j=N}^{\infty}a_j(1+a_j)$.
    If $f \in H^2$ is fixed, then for every $\varepsilon>0$  one can choose $N$ such that
    $3 a^2 \sum\limits_{j=N}^{\infty}a_j(1+a_j) \leq \varepsilon$. Then 
    $g=A_N f \in H^2$ and $T_{N+2n}g \to 0$. 
    Hence, $\|A_{N+2n+1}f\|=\|B_{N+2n+1}A_{N+2n}f\| \leq a \|A_{N+2n}f\| \to 0.$

    So, we have $$f(z)=\sum\limits_{k=1}^{\infty}\sum\limits_{j=1}^{n_k}\frac{c_{k,j}}{1-r_k \bar\zeta_{k,j}z}, \text{ where } \displaystyle c_{k,j}(f)=\int\limits_{I_{k,j}}f_{k-1}(\zeta)dm(\zeta).$$
    Now we have to show that the series $\sum\limits_{k,j}\frac{c_{k,j}}{1-r_k \bar\zeta_{k,j}z}$ converges to $f$.
     To do this, as before, we have to estimate the difference between the intermediate partial sums and the sums over complete circles. The required estimate was established in \cite{BB23} (Section 2, proof of Theorem 1.1).
    \end{proof}

 \begin{Remark}
     We believe that the statement holds in $H^p$ for every $p \in (1,\infty)$, but the above proof crucially relies on the Hilbert structure of $H^2$, namely,
     on the fact that $\|T_n\|_{H^2 \to H^2}=1$. However, $\|T_n\|_{H^p \to H^p}>1$
      if $p \in [1,2)\cup(2,\infty]$. Indeed, the backward shift operator
     $$Bf=(f-f(0))/z$$
     has norm strictly greater than $1$ in any $H^p$ with $p\neq 2$; see \cite{BS24}.
     Let us choose $\varepsilon \in (0, \|B\|_p-1)$, a function $f$ with $\|f\|_{H^p}=1$ such that
      $\|Bf\| >1+\varepsilon$, and put $$g_n(z)=f(z^n).$$
      Note that 
      $\|g_n\|_{H^p}=1$, $\|g_n-g_n(0)\|_{H^p}=
      \|f-f(0)\|_{H^p}$. Then 
      $$\|g_m-g_m(r_n\cdot)\|_{H^p} \ge
       \|g_m-g_m(0)\|_{H^p}-\|g_m(0)-g_m(r_n\cdot)\|_{H^p}$$
       $$=\|f-f(0)\|_{H^p}-\|g_m(0)-g_m(r_n \cdot)\|_{H^p}$$
       $$>1+\varepsilon-\|g_m(0)-g_m(r_n \cdot)\|_{H^p}.$$
       Since $\|g_m(0)-g_m(r_n \cdot)\|_{H^p} \to 0$ as $m \to \infty$, for sufficiently large $M$ one has $$\|g_M-g_M(r_n \cdot)\|_{H^p}>1+\varepsilon/2.$$
       This shows that $\|T_n\|>1$.         
 \end{Remark}

 \subsection{Some generalizations}

Let $g, \, N_k, \, M_k, \, \lambda_{k,l,j}, \, \omega_1, \, \omega_2, \, \omega_3$ be the same as in
Section 4. Using the previous technique, we obtain the following
 \begin{Theorem}
    Let $R_{k, 1} \to 1$ and
$$\lim\limits_{k \to \infty}\frac{\omega_1(R_k) \omega_2(R_k)}{N_k^2}  = 0, \qquad
\omega_1(R_k) \omega_3(R_k) = O(M_k^2), \quad k\to\infty.$$
Then $\{g(\bar\lambda z)\}_{\lambda\in\Lambda}$ is a representing system for $H^2$ and it allows a linear method for the choice of coefficients.
 \end{Theorem}

Let $\gamma=\{\gamma_n\}_{n=0}^{\infty}, \, \gamma_n > 0$,
    $\lim\limits_{n \to \infty} \gamma_n^{1/n}=1$. 
Consider the weighted Hardy space $$\mathscr{H}_\gamma =\left\{f(z)=\sum\limits_{n=0}^{\infty}a_nz^n: \|f\|^2_{\gamma} = \sum\limits_{n=0}^{\infty}|a_n|^2 \gamma_n <\infty \right\}.$$
The space $\mathscr{H}_\gamma$ consists of functions analytic in the unit disk. Moreover, it is an RKHS with
the kernel $$k^{\gamma}_\lambda(z)=\sum\limits_{n=0}^{\infty}\bar{\lambda}^n z^n/ \gamma_n.$$
Weighted Hardy spaces $\mathscr{H}_\gamma$ include most of the classical Hilbert spaces of analytic functions
in the unit disk: the Hardy space $H^2$, the Bergman spaces $A^2_\alpha$,
and the Dirichlet space. If $\inf \gamma_n >0$, then $\mathscr{H}_\gamma$
has no Duffin--Schaeffer frames of normalized reproducing
kernels but has representing systems of reproducing kernels; see \cite{BB23}.

We introduce three characteristics of the kernel: 
$$\theta_1(r)= \sup_{n \ge0}r^{2n} \gamma_n, \quad \theta_2(r)=\sum\limits_{n \ge 1}\frac{n^2 r^{2n}}{\gamma_n},
\quad \theta_3(r)=\sum\limits_{n \ge 0}\frac{r^{2n}}{\gamma_n}.$$

 \begin{Theorem}
   Let $R_{k, 1} \to 1$ and
$$\lim\limits_{k \to \infty}\frac{\theta_1(R_k) \theta_2(R_k)}{N_k^2} = 0, \qquad
\theta_1(R_k) \theta_3(R_k) = O(M_k^2), \, k\to\infty.$$
Then $\{k^{\gamma}_\lambda\}_{\lambda\in\Lambda}$ is a representing system for $\mathscr{H}_\gamma$ and it allows a linear method for the choice of coefficients.
\end{Theorem} 

\begin{Remark}
    By combining the techniques of Sections 4 and 5, one can construct representing systems 
    of dilations of a single function for $\mathscr{H}_\gamma$ that
    allow a linear method for the choice of coefficients.
\end{Remark}

\textbf{Acknowledgement.} The author is grateful to Anton Baranov, Pavel Terekhin, Egor Dobronravov and Nikita Dobronravov for useful discussions.
The work is supported by the Ministry of Science and Higher Education of the Russian Federation (agreement no. 075–15–2025–343).

\printbibliography

\end{document}